\documentclass{article}

\usepackage{arxiv}

\usepackage[utf8]{inputenc} 
\usepackage[T1]{fontenc}    
\usepackage{hyperref}       
\usepackage{url}            
\usepackage{booktabs}   
\usepackage{amsfonts}   
\usepackage{nicefrac}     
\usepackage{microtype}  
\usepackage{graphicx}
\usepackage[numbers,sort&compress]{natbib}
\bibpunct[, ]{[}{]}{,}{n}{,}{,}
\renewcommand\bibfont{\fontsize{10}{12}\selectfont}%
\makeatletter
\def\NAT@def@citea{\def\@citea{\NAT@separator}}
\usepackage{doi}
\usepackage{url}
\usepackage{amsmath,amssymb,amsthm}
\usepackage{interval}
\intervalconfig{soft open fences}
\usepackage{xcolor}
\usepackage{algorithm}
\usepackage{siunitx}
\usepackage{algpseudocode}
\usepackage{hyperref}
\usepackage{tabularx}
\usepackage{comment}

\theoremstyle{plain}
\newtheorem{theorem}{Theorem}[section]

\newtheorem{assumption}[theorem]{Assumption}

\theoremstyle{definition}

\theoremstyle{remark}
\newtheorem{remark}{Remark}

\newcounter{mycount}

\makeatletter
\newcommand{\multiline}[1]{%
  \begin{tabularx}{\dimexpr\linewidth-\ALG@thistlm}[t]{@{}X@{}}
    #1
  \end{tabularx}
}
\makeatother

\title{Piecewise linear value function approximations in nonlinear dynamic scheduling problems with VTOLs}

\author{{Viktoriya Nikitina} \\
	\texttt{viktoriya.nikitina@unibw.de} \\			
	\And
	{Sergejs Rogovs} \\
	\texttt{sergejs.rogovs@unibw.de} \\
	\And
	{Matthias Gerdts} \\
	\texttt{matthias.gerdts@unibw.de} \\
}

\date{Institute of Applied Mathematics and Scientific 		Computing \\
Department of Aerospace Engineering \\
University of the Bundeswehr Munich, Germany \\ \vspace{15pt}
October 29, 2022}
\renewcommand{\shorttitle}{Dynamic Scheduling of VTOLs}

\hypersetup{
pdftitle={A template for the arxiv style},
pdfsubject={q-bio.NC, q-bio.QM},
pdfauthor={David S.~Hippocampus, Elias D.~Striatum},
pdfkeywords={First keyword, Second keyword, More},
}

\begin{document}
\maketitle

\begin{abstract}
Modern vertical take-off and landing vehicles (VTOLs) could significantly affect the future of mobility. The broad range of their application fields encompasses urban air mobility, transportation and logistics as well as reconnaissance and observation missions in a military context. This article presents an efficient numerical framework for dynamic scheduling of multiple VTOLs. It combines a scheduling problem with optimal trajectory planning. The whole setting is formulated as a mixed-integer bilevel optimization problem. At the upper level, VTOLs are scheduled, and their starting times are computed. The solution of the lower level problem involves the computation of a value function and yields optimal trajectories for every aerial vehicle. In order to solve the bilevel problem, it is recast into a single-level one. The resulting mixed-integer nonlinear program (MINLP) is then piecewise linearized and solved numerically by a linear solver based on the Branch-and-Bound algorithm. Numerical results prove the feasibility of the approach proposed in this work.

\end{abstract}


\section{Introduction}
\label{sec:intro}

In the last few years, novel concepts of vertical take-off and landing vehicles (VTOLs) have been developed and successfully tested in practice. Modern approaches are considered to be promising and could considerably affect the future of mobility. The recent technological advances in this field have attracted much attention across different countries and industries. As a result, the global VTOL market has lately seen a significant boom. In Europe, for example, a number of requests have been recently received by the European Union Aviation Safety Agency (EASA) for certifying new types of VTOLs \cite{EASA}.

Thanks to their high maneuverability, comparatively small size and universal deployability (no need for runways), VTOLs are believed to play an important role in the future of several industries. One of them is clearly urban air mobility, where flying taxis could transport passengers from A to B in a much faster way than ground vehicles due to the avoidance of congested streets. Electric VTOLs (eVTOLs) could in addition yield environmental benefits and promote sustainability in the aviation sector. One of the countries blazing a trail in the field of air taxis is Singapore, whose government plans to start the operation of eVTOLs by 2024 at one of the airports \cite{Singapore_Taxi}. Another major area of VTOL application relates to emergency situations, where aerial vehicles could be used for quick medical evacuation of people injured in accidents or for different search-and-rescue scenarios (e.g. after earthquakes, avalanches etc.). Moreover, VTOLs could help to revolutionize the logistics sector. The use of autonomous delivery drones could facilitate the whole shipping process making it more efficient and cost-effective. Multiple key players in the field of electronic commerce have already started to plan the integration of such drones into their transport chain. In a military context, VTOLs could be, for instance, involved in missions aimed at reconnaissance and observation of some targets.

A rapidly increasing number of VTOLs requires their efficient coordination and safe air traffic management, which might be particularly challenging when dealing with limited space. This paper proposes an approach for a small fleet of autonomous VTOLs (typically, up to eight) designed for their simultaneous scheduling and optimal trajectory planning. The task of each VTOL is hereby to reach a static target region without colliding with dynamic obstacles. Provided that the starting position of each aerial vehicle is fixed and given, the overall task is said to be completed as soon as the last VTOL has reached the target set. The whole setting is formulated as a mixed-integer bilevel optimization problem. At the upper level, the scheduling of multiple autonomous VTOLs can be viewed as a job-shop problem (JSP) with the aim of finding an optimal sequence and the corresponding starting times. At the lower level, the task is to determine the shortest path from every aerial vehicle to the target region. The solution procedure involves computation of a value function and yields an optimal trajectory for every VTOL provided that the starting times are given. Clearly, the optimization problems at both levels are strongly coupled with each other and have to be therefore considered as a whole.

To date, different efficient approaches have been developed to solve the optimization problems at the upper and lower level if they are considered separately. The authors of \cite{Parzani2017} and \cite{Coupechoux2021} deal, for instance, with optimal trajectory planning for aircraft using the Hamilton-Jacobi-Bellman (HJB) equation. \cite{Rogovs2022} and \cite{Chandran2007} in turn present exact methods for coordination of multiple landing aircraft. However, the class of \textit{coupled} scheduling and optimal control problems has not been extensively researched in the literature to the best of the authors' knowledge. In a context of autonomous vehicles, \cite{Mahdavi2018} and \cite{Gerdts2022} propose interesting methods for optimization problems with a bilevel structure. Approaches presented in \cite{Palagachev2017} and \cite{Zanlongo2018} are designed to schedule multiple robots and provide an optimal path for each of them. In \cite{Sama2017}, a bilevel mixed-integer optimization problem is investigated with the aim of scheduling aircraft and planning their optimal trajectory in a busy terminal area.

The remainder of this article is structured as follows. Section \ref{sec:opt_con_HJB} briefly introduces the value function of an optimal control problem, which solves the Hamilton-Jacobi-Bellman (HJB) equation. This theoretical background originates from dynamic programming and is required in succeeding sections. The mathematical model used for dynamic scheduling of VTOLs is formulated in Section \ref{sec:math_problem}. Its solution procedure is then carefully described in Section \ref{sec:sol_procedure}. Herein, one of the steps includes an elegant technique used in this article for the piecewise linearization of the final MINLP. The solution of the resulting linearized problem can be obtained numerically using the robust and efficient GUROBI solver. Section \ref{sec:numerics} is mainly devoted to the presentation of the results achieved by the approach presented in this work. Finally, some concluding remarks and possible directions for further research are given in Section \ref{sec:conclusion}.

\section{Optimal Control and Value Function}
\label{sec:opt_con_HJB}

The main goal of this section is to give a short overview of the theory used in Sections \ref{sec:math_problem} and \ref{sec:sol_procedure}. The concepts briefly outlined below deal with the Hamilton-Jacobi-Bellman (HJB) equation and the value function approach. While Section \ref{subsec:HJB_cont} discusses the continuous time problem, Section \ref{subsec:HJB_discr} presents formulas related to the discrete case. Note that the derivation of formulas is omitted here since it goes beyond the scope of this article. More detailed information can be, for example, found in \citep{Bardi1997, Bellman1957}.

\subsection{Continuous Time}
\label{subsec:HJB_cont}

For $\lambda \geq 0$ being a given constant discount rate, $U \subset \mathbb{R}^{n_u}$ a control space (which is typically closed and bounded), $x \in W^{1,\infty}\left(\interval[open right]{0}{\infty}, \mathbb{R}^{n_x}\right)$ a state function and $u \in L^{\infty}\left(\interval[open right]{0}{\infty},U\right)$ a control function, consider the following continuous optimal control problem (OCP) on an infinite time horizon:

\begin{problem}{Continuous OCP}
	Minimize
	\begin{equation}
	\label{prob:opt_con_cts}
		J(x_0, u) := \int_{0}^{\infty} {\exp(-\lambda t) C(x(t), u(t)) dt}
	\end{equation}
	subject to the constraints
	\begin{equation*}
		\begin{array}{rcl}
			x^{\prime}(t) &= &F(x(t), u(t)),\\
			x(0) &= &x_0, \\
			u(t) &\in &U.
		\end{array}
	\end{equation*}
\end{problem}
Herein, $x_0$ is a given initial state, $F: \mathbb{R}^{n_x} \times U \to \mathbb{R}^{n_x}$ is a continuous function describing the state dynamics, and $C: \mathbb{R}^{n_x} \times U \to \mathbb{R}$ is a continuous running cost function.

Let the value function $V$ of the optimal control problem \eqref{prob:opt_con_cts} be defined as 
\begin{equation*}
\label{def_val_fct_cts}
	V(x_0) := \inf\{J(x_0, u) \mid u(t) \in U \ \text{for almost every} \ t \in \mathbb{R}_{\geq 0}\}.
\end{equation*}
It can be shown \cite{Bardi1997} that $V$ solves the following equation
\begin{equation*}
\label{eq:HJB_cts}
	\lambda V(x) + \sup_{u\in U}\{ -\nabla V(x)^T F(x,u) - C(x,u) \} = 0,
\end{equation*}
which is called the Hamilton-Jacobi-Bellman (HJB) equation \cite{Bellman1957}.

For a given $x \in \mathbb{R}^{n_x}$, the feedback law can be then computed as 
\begin{equation*}
	u(x) = \underset{u \in U}{\arg \min} \{ \nabla V(x)^T F(x,u) + C(x,u)\}, 
\end{equation*}
whenever this minimum exists. Note that the gradient of $V$ is denoted by $\nabla V(x)$ for $V$ supposed to be sufficiently smooth.

\subsection{Discrete Time}
\label{subsec:HJB_discr}
Consider the following discrete optimal control problem (OCP) with a constant discount rate $0 < \beta \leq 1$:

\begin{problem}{Discrete OCP}
	Minimize
	\begin{equation}
		\label{prob:opt_con_discr}
		j(x_0, \{u_n\}_{n \in \mathbb{N}_0}) := \sum_{n=0}^{\infty} {\beta^n c(x_n, u_n)}
	\end{equation}
	subject to the constraints
	\begin{equation*}
		\begin{array}{ll}
			\begin{array}{rcl}
				x_{n+1} &= &f(x_n, u_n), \\
				u_n &\in& U
			\end{array} & \text{for} \ n = 0, 1, \dots
		\end{array}
	\end{equation*}
\end{problem}
Herein, $f: \mathbb{R}^{n_x} \times U \to \mathbb{R}^{n_x}$ is a continuous function describing the discrete state dynamics, and $c: \mathbb{R}^{n_x} \times U \to \mathbb{R}$ is a continuous running cost function.

A detailed look at the objective function of problem \eqref{prob:opt_con_discr} reveals that it is constructed in such a way that the closer the summands to the initial state $x_0$ are, the stronger their influence on $j$ is due to the weighting parameters $\beta^n, \ n \in \mathbb{N}_0$.

The value function $V_h$ of the discrete optimal control problem \eqref{prob:opt_con_discr} is defined as
\begin{equation}
\label{def:val_fct_discr}
	V_h(x_0) := \inf\{ j(x_0, \{u_n\}_{n \in \mathbb{N}_0}) \mid u_n \in U \ \text{for} \ n \in \mathbb{N}_0 \}
\end{equation}
and represents the greatest lower bound of the total cost for all trajectories, whose starting points are $x_0$ \cite{Schreiber2016}. It can be shown \cite{Bellman1957} that $V_h$ solves the following fixed-point equation for $V_h(x)$
\begin{equation}
	\label{eq:HJB_discr}
	V_h(x) = \inf_{u \in U}\{ c(x,u) + \beta V_h(f(x,u)) \},
\end{equation}
which is a discrete version of the HJB equation. The derivation of the equation exploits the definition of $V_h$ and involves, inter alia, the insertion of equation \eqref{prob:opt_con_discr} into equation \eqref{def:val_fct_discr}. For reasons of space, it is omitted here. The interested reader is instead referred to \cite{Schreiber2016} for more details.

Equation \eqref{eq:HJB_discr} can be solved using the following fixed-point iteration for an initial guess $V_h^0$
\begin{equation}
\label{iter_HJB_discr}
	V_h^{m+1}(x) = \inf_{u \in U}\{ c(x,u) + \beta V_h^m(f(x,u)) \} \quad \text{for} \ m \in \mathbb{N}_0
\end{equation}
until the sequence of functions $\{V_h^m\}_{m \in \mathbb{N}_0}$ converges \cite{Bardi1997}. 

For $\beta \in \interval[open]{0}{1}$, the operator
\[
	T[V](x) := \inf_{u \in U} \{ c(x,u) + \beta V(f(x,u))\}
\] 
is a $\beta$-contraction. Thus, equation \eqref{iter_HJB_discr} possesses a unique fixed point according to the Banach fixed-point theorem. However, it cannot be guaranteed that the discrete HJB equation has a unique fixed point for $\beta = 1$. An elegant way of handling this issue is the so-called Kru\u{z}kov transform \cite{Kruzhkov1975} discussed later in Section \ref{sec:first_step}.

\begin{remark}
\label{rem:trafo}
	Note that the explicit dependency on time $t$ of a non-autonomous OCP can be eliminated by introducing a new state variable
\[ 
	\tau(t) := t. \\
\]
In the continuous case, the new state vector corresponds then to
\[
	\tilde{x}(t) = \begin{bmatrix} x(t) \\ \tau(t) \end{bmatrix},
\]
which solves the following initial value problem (IVP):
\begin{align*}
	\tilde{x}'(t) &= \begin{bmatrix} F( \tau(t), x(t), u(t)) \\ 1 \end{bmatrix}, \\
	\tilde{x}(0) &= \begin{bmatrix} x_0 \\ 0 \end{bmatrix}.
\end{align*}
Similarly, in the discrete case it holds
\begin{align*}
	\tilde{x}_{n+1} = \begin{bmatrix} x_{n+1} \\ \tau_{n+1}  \end{bmatrix}  = \begin{bmatrix} f(x_n,u_n) \\ \tau_n + h \end{bmatrix} \quad \text{for} \ n \in \mathbb{N}_0,
\end{align*}
where $h > 0$ is the step size of the time domain grid assumed to be equidistant. Therefore, a value function $V_h(t,x(t))$ explicitly depending on $t$ can be transformed into a function $V_h(\tilde{x}(t))$ depending solely on the state variables.
\end{remark}

\subsection{Relation between Continuous and Discrete Optimal Control Problems}
\label{subsec:relation_con_discr}

In order to show a relation between the continuous OCP and its discrete counterpart, the continuous problem is discretized on the grid $t_n = nh, \ n \in \mathbb{N}_0,$ with a step size $h > 0, \ x_n := x(t_n)$ and $u_n := u(t_n)$. This implies:
\begin{equation*}
	J(x_0, u) = h \sum_{n = 0}^{\infty} {\exp(-\lambda t_n) C(x_n, u_n)} = \sum_{n = 0}^{\infty} {{\exp(-\lambda h)}^n h C(x_n, u_n)}.
\end{equation*}
Hence, the relation holds with $\beta = \exp(-\lambda h)$ and $c(x,u) = h C(x,u)$.

To obtain a relation between the functions $f$ and $F$, the Euler method is exploited for dealing with the ODE describing the state dynamics. This results in
\begin{equation*}
	F(x_n, u_n) = x^{\prime}(t_n) \approx \frac{x_{n+1} - x_n}{h} = \frac{f(x_n, u_n) - x_n}{h}.
\end{equation*}
Thus, it can be recognized that $f(x,u) = x + hF(x,u)$.

\section{Problem Formulation}
\label{sec:math_problem}

This section is devoted to the exact formulation of the problem. The overall goal is to determine an optimal schedule for a small fleet of autonomous VTOLs and to find the respective trajectories. Hereby, their dynamics and the condition that while one VTOL executes its operation, all others have to idle are taken into account. In our case, the mission for all VTOLs is to reach a static target from some starting points while avoiding moving obstacles such that the total mission time, i.e. the sum of the respective flight durations, is minimal. Note that flight durations can be penalized.

This setting can be formulated as a mixed-integer bilevel optimization problem, where the upper level problem deals with the scheduling of autonomous VTOLs, while the lower level one delivers optimal trajectories. In what follows, some notation used in this paper is introduced. Afterwards, the underlying problem's mathematical formulation is discussed in detail.

\begin{figure}[htbp!]
\begin{center}
\includegraphics[scale=1]{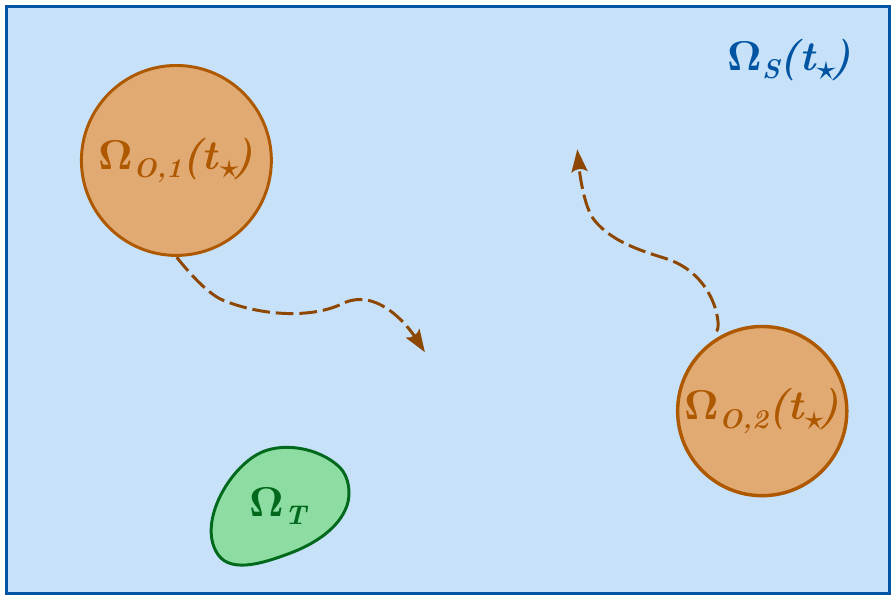}
\end{center}
\caption{Example of $\Omega(t_\star),\ t_\star \in \mathbb{R}_{\geq 0}$. In green: static target set $\Omega_T$, in orange: dynamic obstacles with $\Omega_O(t_\star) = \Omega_{O,1}(t_\star) \cup \Omega_{O,2}(t_\star)$, in blue: time-dependent set $\Omega_S(t_\star)$, across which VTOLs can move.}
\label{fig:omega_set}
\end{figure}

Throughout this work, the number of autonomous VTOLs is denoted by $N \in \mathbb{N}$ with the corresponding index set $\mathcal{N}:= \{1, \dots, N\}$. Moreover, $P \in \mathbb{N}_0$ and $\mathcal{P}:= \{1, \dots, P\}$ indicate the number of (moving) obstacles, whose dynamics is supposed to be known, and its index set, respectively. In addition, the fixed starting positions of all aerial vehicles are assumed to be given. The overall goal is to schedule VTOLs and to simultaneously compute their optimal trajectories in some time-dependent domain $\Omega(t) \subset \mathbb{R}^{n_x}$ with $ 
t \in \mathbb{R}_{\geq 0}$ and $n \in \{2,3\}$ such that all of them reach a predefined static target set $\Omega_T \subset \Omega(t), \ t \in \mathbb{R}_{\geq 0},$ while avoiding possible collisions with any of the obstacles. Note that the latter form the time-dependent set $\Omega_O(t)$ with $\Omega_O(t) = \underset{i \in \mathcal{P}}{\cup} \Omega_{O,i}(t) \subset \Omega(t), \ t \in \mathbb{R}_{\geq 0}$. The set, across which VTOLs are allowed to move at any time point, is denoted by $\Omega_S(t)$. The whole domain $\Omega(t)$ can be thus represented as $\Omega(t) = \Omega_T \cup \Omega_O(t) \cup \Omega_S(t), \ t \in \mathbb{R}_{\geq 0}$. An illustrative example of $\Omega(t_\star)$ with two obstacles is given in Figure \ref{fig:omega_set} for a specific time point $t_\star \in \mathbb{R}_{\geq 0}$. 

\vspace{5pt}
In addition, the following assumptions must apply within the framework of this article. We briefly comment each of them.

\begin{assumption}
\label{ass:fin_time}
	The target region $\Omega_T$ can be reached by every VTOL in finite time.
\end{assumption}
This assumption is introduced in order to exclude situations where, for instance, the target set is unreachable because of insurmountable obstacles.

\begin{assumption}
\label{ass:const_velocity}
	The velocity of all VTOLs is assumed to be constant.
\end{assumption}
This assumption is legitimate, since a VTOL moves with a constant velocity most of the flight. It allows to consider a simplified version of the underlying OCPs, see Problem \eqref{prob:OCP}.

\begin{assumption}
\label{ass:t*}
	All obstacles stop their motion at some time point $t^{\star}$ bounded from below by the time point, at which the last VTOL in the sequence reaches $\Omega_T$.
\end{assumption}
This is a reasonable assumption, since the overall mission is said to be fulfilled after the last aerial vehicle has reached the target set. It is essential for our computational approach, see Section \ref{sec:first_step}.

\vspace{5pt}
We are now at a position to formulate the whole problem with a bilevel structure. As already mentioned before, the VTOLs in our consideration are allowed to execute their respective operations only one by one. Hence, the scheduling part of the problem can be formulated as the so-called job-shop problem (JSP) \cite{Blazewicz1996}.

For reasons of comprehensiveness, the bilevel problem with a coupling between the two levels is presented  as an upper and lower level problem and can be formulated as 
\begin{problem}{Upper Level: Scheduling of VTOLs}
\begin{subequations}\label{prob:SP}
	\begin{equation}
		\min_{t_i, w_{ij}} \sum_{i=1}^{N} \left[ t_i + \alpha D_i(t_i) \right]
	\end{equation}
\noindent
	subject to
	\begin{alignat}{4}
 		\label{constr:sch1}
		&& t_i + D_i(t_i) - t_j \quad && \leq \quad &&& M(1 - w_{ij}), \\
		\label{constr:sch2}
		&& t_j + D_j(t_j) - t_i \quad && \leq \quad &&& M w_{ij}, \\
		\label{constr:t}
		&& t_i \quad && \in \quad &&& [t_{i,min},t_{i,max}], \\
		\label{constr:wij}
		&& w_{ij} \quad && \in \quad &&& \{0,1\}, \ i,j \in \mathcal{N}, i<j.
	\end{alignat}
\end{subequations}
\end{problem}
Herein, $t_i$ denotes the starting time of VTOL $i, \ i \in \mathcal{N},$ $\alpha$ is a flight duration penalization term, $M$ is some sufficiently large constant, which guarantees that only one of constraints \eqref{constr:sch1}, \eqref{constr:sch2} is active at the same time. The binary variables  $w_{i,j}, \ i,j \in \mathcal{N},$ have the following meaning:
\[
  w_{ij} =
    \begin{cases}
        1 & \textnormal{if VTOL $i$ starts before VTOL $j$},\\
        0 & \textnormal{otherwise}.
    \end{cases}     
\]
Moreover, $D_i( t_i), \ i \in \mathcal{N},$ is a nonlinear function and denotes the flight duration of the corresponding VTOL if it starts at time $t_i$, which in turn is a continuous variable bounded by $t_{i,min}$ and $t_{i,max}$ from below and above, respectively.
Furthermore, the following relation holds due to Assumption \ref{ass:const_velocity}
\[
D_i(t_i) = d(t_i)/v_i, \quad i \in \mathcal{N},
\]
where $v_i$ is the average velocity of VTOL $i$. For simplicity reasons, we assume that $v_i = v,\ \forall i \in \mathcal{N}$. The numerator in the previous equality is a parametric function, which denotes the length of the shortest path to the target at time $t_i,\ i \in \mathcal{N}$, and can be computed from the following parametric optimal control problem:

\begin{problem}{Lower Level: Shortest Path}
\begin{subequations}\label{prob:OCP}
	\begin{equation}
			d(t_i) := \min_{x_i, u_i} 
			  \int_{t_i}^{\infty} {\chi_{\Omega_S(t)}(x_i(t)) dt}
	\end{equation}
subject to 
	\begin{alignat}{4}
		\label{constr:OCP_1}
		& x_i^{\prime}(t) && = \  && u_i(t) \quad && \textit{for a.e.} \ t \in 		\interval[open right]{t_i}{\infty},  \\
		\label{constr:OCP_2}
		& x_i(t)  && \in \ && \Omega_S(t)\cup \Omega_T \quad && \textit{for every} \ t \in 			\interval[open right]{t_i}{\infty}, \\
		\label{constr:OCP_3}
		& u_i(t) && \in \ && U_i(x_i) = \begin{cases} \{ u \in \mathbb{R}^{n_u} \mid || u || = 1\}, \ x_i \in \Omega_S(t) \\ \{0 \in \mathbb{R}^{n_u} \}, \ x_i \in \Omega_T \end{cases} \quad && \textit{for a.e.} \ t \in \interval[open right]{t_i}{\infty}, \\
		\label{constr:OCP_4}
		& x_i(t_i)  && = \ &&  \overline{x}_i, \quad && 
\end{alignat}
\end{subequations}
\end{problem}
where $\overline{x}_i$ is the initial position of VTOL $i$, and $\chi_A$ is the indicator function defined as
\begin{equation*}
	\chi_A(x) =
	\begin{cases}
		1 & \text{if} \ x \in A, \\
		0 & \text{otherwise}.
	\end{cases}
\end{equation*}

Note that in case that $n_u=2$, the control sets $U_i(x_i)$ for $x_i \in \Omega_S(t),\ i \in \mathcal{N},$ can be represented in the following way using polar coordinates:
\begin{equation*}
	U_i(x_i) = \left\{ \left( \
	\begin{matrix}
		u_1 \\
		u_2
	\end{matrix} \ \right)
	\ \Big| \ u_1 = \cos{\alpha_i},\ u_2 = \sin{\alpha_i},\ \alpha_i \in \interval[open right]{0}{2\pi} \right\}.
\end{equation*}

Similarly, the spherical coordinate system can be exploited for describing $U_i(x_i)$ for $x_i \in \Omega_S(t),\ i \in \mathcal{N},$ if $n_u = 3$: 
\begin{equation*}
	U_i(x_i) = \left\{ \left( \
	\begin{matrix}
		u_1 \\
		u_2 \\
		u_3
	\end{matrix} \ \right)
	\ \bigg| \ u_1 = \cos{\alpha_i} \sin{\beta_i},\ u_2 = \sin{\alpha_i} \sin{\beta_i},\ u_3 = \cos{\beta_i}\, \alpha_i \in \interval[open right]{0}{2\pi},\ \beta_i \in \interval{0}{\pi} \right\}.
\end{equation*}

By setting 
\begin{align*}
	F(x,u) &= u, \\
	C(x,u) &= \chi_{\Omega_S(t)}(x), \\
	\lambda &= 0,
\end{align*}
it becomes clear that this problem fits into the framework of OCP classes, see \eqref{subsec:HJB_cont}. Hence, it can be solved with the value function approach. This is very beneficial in our case, since only one value function has to be computed for all problems due to their structure. In this case, it has to be computed for $t \in \mathbb{R}_{t \geq t_{min}},\ t_{min} = \min_{i \in \mathcal{N}} t_{i,min}$. Without loss of generality, we will hereinafter always assume $t_{min} = 0$.

Obviously, problem \eqref{prob:OCP} is time-dependent because the underlying domain depends on time. In order to fit this problem into the framework of OCP classes, steps from Remark \ref{rem:trafo} have to be exploited. However, they are omitted here due to readability reasons and appear only in the discrete case, see Section \ref{sec:first_step}. It can be shown that solving  problem \eqref{prob:OCP} is equivalent to solving the Eikonal equation given by:
\[
\label{eq:Eikonal}
\begin{split}
	 ||\nabla V_h(x_i)|| &= 1 \quad \text{for} \ x_i \in \Omega_S(t),\ t \in \mathbb{R}_{t \geq 0}, \\
	V_h(x_i) &= 0 \quad \text{for} \ x \in \partial \Omega_T, t \in \mathbb{R}_{t \geq 0}.
\end{split}
\]
Therefore, the solution of this problem yields the shortest paths from $\overline{x}_i, \ i \in \mathcal{N},$ to the target set $\Omega_T$.

\section{Solution Procedure}
\label{sec:sol_procedure}
This section highlights and explains main steps required to solve the bilevel mixed-integer optimization problem \eqref{prob:SP} - \eqref{prob:OCP}. The solution procedure can be divided into three main steps indicated below:
\begin{enumerate}
	\item Solve the lower level problems corresponding to the parametric optimal control problems \eqref{prob:OCP}. Its output, as discussed above, is a single value function $V_h(t,x)$ depending on both the starting time $t$ and the starting position $x$.
	\item Insert all given VTOLs' starting positions $\overline{x}_i, \ i \in \mathcal{N},$ into the value function $V_h$ obtained in step (1) and set $d_i(t_i) := V_h(t_i, \overline{x}_i)$ to obtain time-dependent shortest-path functions. The bilevel problem \eqref{prob:SP} - \eqref{prob:OCP} is hence recast into a single-level MINLP. It is then piecewise linearized and solved by a robust and efficient MILP-solver. A solution to that problem is an optimal sequence of VTOLs provided by the optimal starting times $t^{\ast}_i, \ i \in \mathcal{N}$.
	\item Compute the optimal trajectory of every VTOL via the forward recursion.
\end{enumerate}

\subsection{First Step: Solving the Optimal Control Problem}
\label{sec:first_step}
\begin{figure}[htbp!]
\begin{center}
\includegraphics[scale=1]{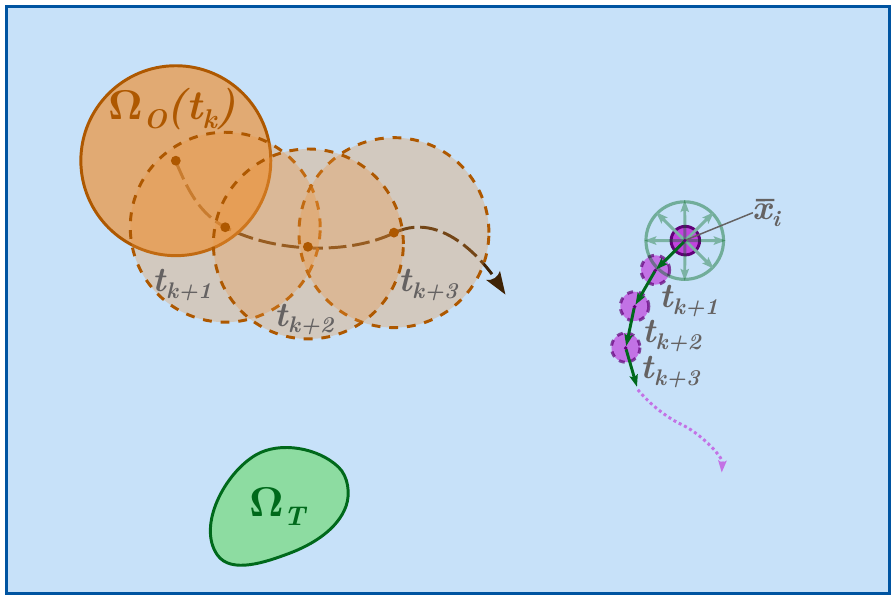}
\end{center}
\caption{Approach for the shortest path computation. In green: static target set $\Omega_T$, in orange: a dynamic obstacle $\Omega_O(t)$, in purple: a moving VTOL}
\label{fig:distance}
\end{figure}

In the first step of the solution procedure, the parametric optimal control problems \eqref{prob:OCP} need to be solved. Before we go deeper into detail, let us first explain the idea behind the formulation of \eqref{prob:OCP}. Starting at a point $\overline{x}_i$, a circle with a radius of one and a midpoint at $\overline{x}_i$ is considered. The aim is to find that point lying on the boundary of this circle, which has the minimum distance to the target set, avoiding a collision with dynamic obstacles. The distance between $\overline{x}_i$ and the newly found point is thus equal to one. In the next step, a new circle with its midpoint lying in this new point is considered. The whole procedure is repeated until the last circle crosses the target area. In the end, the length of the shortest path from point $\overline{x}_i$ to the target set $\Omega_T$ can be computed as a sum of ones. Figure \ref{fig:distance} visualizes the first steps during the computation procedure according to this idea.

To obtain a solution of problem \eqref{prob:OCP}, the value function approach discussed in Section \ref{subsec:HJB_discr} is exploited. Although it is known to suffer from the "curse of dimensions", in the case of the shortest-path problem the dimensionality is characterized only by the time and space variables. Thus, the problem possesses minimum possible dimensionality. Moreover, since the value function is only required to obtain the optimal starting times, simple VTOL models can be used to determine them. Afterwards, more complex models can be exploited in order to compute more precise trajectories. Note that the corresponding OCPs will be no longer parametric in this case and can be therefore solved with standard techniques.

To determine the setting of the discrete problem, we use the relation between it and its continuous counterpart explained in Section \ref{subsec:relation_con_discr} and the argument from Remark \ref{rem:trafo}. As a result, we obtain: 
\begin{equation*}
	\begin{array}{rcl}
		f(x,u) &= & x + hu, \\
		c(x,u) &= & \begin{cases} h & \text{if} \ x \in \Omega_S(t_k), \\ 0 & \text{otherwise}, \end{cases} \\
		\beta \quad &= & 1, \\
		t_{k+1} &= & t_k + h, \quad k \in \mathbb{N}_0.
	\end{array}
\end{equation*}

For the value function, it holds
\begin{equation}\label{eq:v_start}
	V_h(t_k,x) = 
	\begin{cases} \infty & \text{for} \ x \in \left( \mathbb{R}^{n_x} \setminus \Omega \right) \cup \Omega_O(t_k), \\
	0 & \text{for} \ x \in \Omega_T
	\end{cases} \quad \text{for} \ k \in \mathbb{N}_0.
\end{equation}

The value function can be computed iteratively using formula \eqref{iter_HJB_discr} indicated in Section \ref{sec:opt_con_HJB}.
However, this might cause several difficulties. From a theoretical point of view, it is not clear a priori whether equation \eqref{eq:HJB_discr} possesses a unique solution for $\beta = 1$. From a numerical point of view, the fact that $V_h$ can attain infinite values poses another problem, since they might be propagated \cite{Schreiber2016}. An effective remedy is provided by the so-called Kru\u{z}kov transform \cite{Kruzhkov1975, Schreiber2016} defined as:
\begin{equation}\label{eq:kruzkov} 
	v_h(t_k, x):= \exp(-V_h(t_k,x)), \quad  k \in \mathbb{N}_0. 
\end{equation}
Applying \eqref{eq:kruzkov} to	\eqref{eq:HJB_discr} and taking into account \eqref{eq:v_start} yields
\begin{equation*}
	v_h(t_k,x) = 
	\begin{cases} 
		0 & \text{for} \ x \in \left( \mathbb{R}^{n_x} \setminus \Omega \right) \cup \Omega_O(t_k), \\
		\underset{u \in U}{\sup}\{ \exp(-c(x,u)) v_h(t_{k+1},f(x,u)) \} & \text{for} \ x \in \Omega_S(t_k), \\
		1  & \text{for} \ x \in \Omega_T
	\end{cases}
\end{equation*}
for $k \in \mathbb{N}_0$.

Similarly to \eqref{iter_HJB_discr}, the transformed value function can be computed iteratively
\begin{equation*}
	v_h^{m+1}(t_k,x) := \sup_{u \in U}\{ \exp(-c(x,u)) v_h^m(t_{k+1},f(x,u)) \} \quad \text{for} \ k,m \in \mathbb{N}_0
\end{equation*}
with the following initial guess
\begin{equation*}
	v_h^0(t_k,x) = 
	\begin{cases}
		1 & \text{if} \ x \in \Omega_T, \\
		0 & \text{otherwise}
	\end{cases} \quad \text{for} \ k \in \mathbb{N}_0.
\end{equation*}

Provided that $c(x,u)$ is positive, the sequence of functions $\{ v_h^m \}_{m \in \mathbb{N}_0}$ will converge to some function $\hat{v}_h$. The original value function can be then obtained using the following backtransformation:
\begin{equation*}
	V_h(t_k,x)= -\ln\big(\hat{v}_h(t_k,x)\big), \quad k \in \mathbb{N}_0.
\end{equation*}

One of problems, which arises during the numerical computation of the value function $V_h$, relates to the question how the infinite dimensional problems can be handled in practice. To overcome this issue, the following two cases are considered:

\begin{enumerate}
	\item All obstacles are supposed to move periodically with some global period $T$. In this case, the value function can be computed in the following way:
		\begin{align*}
			& v_h^{m+1}(t_k, x) &&=\sup_{u \in U}\{ \exp(-c(x,u)) v_h^m(t_{k+1},f(x,u)) \} & \text{for} \ k = 0, \dots, K-1, \\
			& v_h^{m+1}(t_{K},x) &&= \sup_{u \in U}\{ \exp(-c(x,u)) v_h^m(t_{0},f(x,u)) \}, &	
		\end{align*}
	where $Kh = T$.	
	\item Due to Assumption \ref{ass:t*}, the underlying domain can be viewed as static after some time point $t^{\star}$, hence 
\begin{align*}
	& v_h^{m+1}(t_k,x) &&= \sup_{u \in U}\{ \exp(-c(x,u)) v_h^m(t_{k+1},f(x,u)) \} & \text{for} \ k = 0, \dots, K-1,
	\\
	& v_h^{m+1}(t_{K},x) &&= \sup_{u \in U}\{ \exp(-c(x,u)) v_h^m(t_{K},f(x,u)) \}, &	
\end{align*}
where $Kh > t^{\star}$.

\end{enumerate}

In numerical computations, the space domain is discretized with some equidistant gridpoints. To obtain values of the value function in non-grid points, piecewise linearization can be exploited. 

\subsection{Second Step: Solving the Scheduling Problem}

Using the value function obtained in the previous step, the duration functions of each VTOL can be now defined as
\[
	d_i(t_i):= V_h(t_i, \overline{x}_i), \quad i \in \mathcal{N}.
\]
Recall that $\overline{x}_i$ is the initial starting position of VTOL $i, \  i \in \mathcal{N}$. Hence, one ends up with MINLP \eqref{prob:SP}.
\newline
In general, MINLPs are very hard to solve, especially for very nonlinear problems. To overcome this issue, this work proposes an approximation technique with piecewise linear functions. It eventually yields an MILP and allows to track the linearzation error.

To this end, define a closed interval in $\mathbb{R}$ by
\[
	I_{a}^{b} := \interval{a}{b}.
\]
For time intervals $I_{t_{i,min}}^{t_{i,max}},\ i \in \mathcal{N},$ consider their discretizations
\[
\mathcal{G}_i := \{t_{i,1}, t_{i,2}, \dots, t_{i,K_i} \},
\] 
where $t_{i,1} = t_{i,min}, \ t_{i,K_i} = t_{i,max}$, and $K_i$ denotes the number of discretization points in the grid $\mathcal{G}_i, \ i \in \mathcal{N}$.

Next, consider a nonlinear function $g: \interval{t_a}{t_b} \to \mathbb{R}$ and its linear  approximation $\psi(t)$ satisfying the conditions $g(t_a) = \psi(t_a)$ and $g(t_b) = \psi(t_b)$. The maximum \textit{underestimator} error is then defined as
\[
e_{u}(g,I_{t_a}^{t_b}) := \max_{t \in I_{t_a}^{t_b}}  g(t) -  \psi(t),
\]
whereas the maximum \textit{overestimator} error is indicated by
\[
e_{o}(g,I_{t_a}^{t_b}) := \max_{t \in I_{t_a}^{t_b}} \psi(t) - g(t).
\]
For more details, the interested reader is referred to \cite{Gerdts2022} and \cite{Burlacu2019_PhD}.

\begin{figure}[htbp!]
\begin{center}
\includegraphics[scale=1]{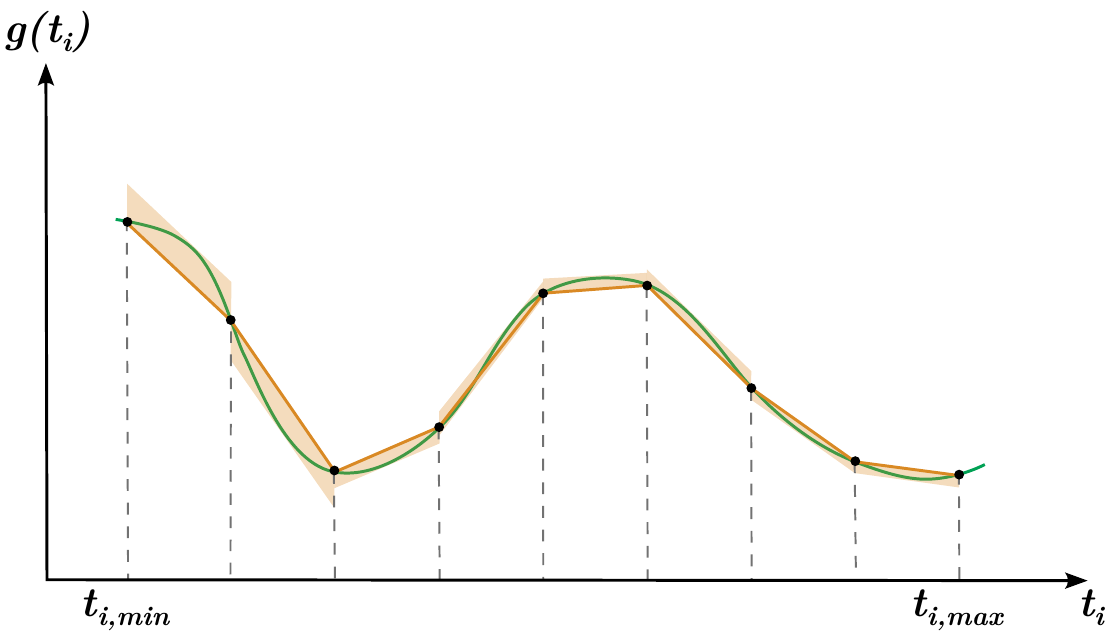}
\end{center}
\caption{Piecewise linearization with envelope functions}
\label{fig:pw_lin_fun}
\end{figure}

Figure \ref{fig:pw_lin_fun} shows a piecewise linearization of a nonlinear function $g(t_i)$, whereby the so-called \textit{envelopes} around the piecewise linear function are highlighted in orange.

In order to obtain a fully linear problem, which can be directly handled by a linear solver, two sets of auxiliary variables are additionally introduced. Namely, binary variables $z_{i,k}$ and real variables $\lambda_{i,k}$ for $i \in \{ 1, \dots, N \}$ with $k \in \{ 1, \dots, K_i-1 \}$ and $k \in \{ 1, \dots, K_i \}$, respectively. The integer variables $z_{i,k}$ are called \textit{pointers} and indicate that subinterval, within which the optimal value of $t_i$ lies. As to the real variables $\lambda_{i,k}$, they represent the coefficients in the optimal subinterval required to build a convex combination there. Obviously, $t_i$ attains its optimal value on exactly one subinterval, wherefore there is only one non-zero pointer. 

These arguments lead to the following constraints:

\begin{align}
	\begin{split}
	\label{eq:PW_contraints_z}
	& \sum_{k=1}^{K_i-1}z_{i,k} = 1,\\
	& z_{i,k} \in  \left \lbrace 0,1 \right \rbrace, \quad k \in \{ 1, \dots, K_i-1 \}
	\end{split}
\end{align}
and
\begin{align}
	\begin{split}
	\label{eq:PW_contraints_lam}
		\sum_{k=1}^{K_i}\lambda_{i,k} = 1.
	\end{split}
\end{align}
In addition, the following inequalities have to be incorporated into the model to guarantee that two consecutive lambdas are unequal to zero:
\begin{align}
	\begin{split}
	\label{eq:PW_contraints_z_lam}
		\lambda_{i,1} & \leq  z_{i,1}, \\
		\lambda_{i,k} & \leq  z_{i,k-1} + z_{i,k}, \quad k \in \{ 2, \dots, K_i-1 \},\\
		\lambda_{i,K_i} & \leq  z_{i,K_i-1}.\\
	\end{split}
\end{align}
Thus, the following additional constraints result from the preceding considerations regarding the piecewise linearization technique for $i \in \{ 1, \dots, N\}$:
\begin{subequations}
	\label{eq:PW_all}
	\begin{alignat}{4}
	t_i &= \sum_{k=1}^{K_i} \ \lambda_{i,k} \ t_{i,k}, \label{eq:PW_ti} \\
	D_{i,pw} &:= \sum_{k=1}^{K_i} \ \lambda_{i,k} \ D_i(t_{i,k}) + e_{i}, \label{eq:PW_D} \\
	e_{i} &\leq \sum_{k=1}^{K_i-1} e_{u}\left(D_i,I_{t_{i,k}}^{t_{i,k+1}}\right) z_{i,k}, \label{eq:under_D}\\
	e_{i} &\geq - \sum_{k=1}^{K_i-1} e_{o}\left(D_i,I_{t_{i,k}}^{t_{i,k+1}} \right) z_{i,k}. \label{eq:over_D}
	\end{alignat}
\end{subequations}
Herein, \eqref{eq:PW_ti} represents the starting time variable $t_i$ as a linear combination of discretization points, \eqref{eq:PW_D} corresponds to the linearizations of the nonlinearities $D_i(\cdot)$ with the corresponding envelope terms. Those are bounded from above and below by the maximum under- and overestimators \eqref{eq:under_D}-\eqref{eq:over_D}, respectively.

With this in mind, the mixed-integer linear relaxation of problem  \eqref{prob:SP} can be finally formulated as 
\begin{problem}{MILP}
	\begin{equation}
		\tag{$\Pi$}
		\label{prob:MILP}
		\begin{array}{l}
			\displaystyle\min_{t_i,w_{ij}, e_{i} }\ \sum_{i=1}^{N} \left[ t_i+ \alpha D_{i,pw} \right],\\
			\textnormal{s.t.} \quad 
			\eqref{constr:sch1} - \eqref{constr:wij},\eqref{eq:PW_contraints_z} -\eqref{eq:PW_all}.
		\end{array}
	\end{equation}
\end{problem}
Herein, all continuous variables in constraints \eqref{constr:sch1} - \eqref{constr:t} are replaced by their discrete counterparts.

The formulation of the adaptive numerical algorithm used for the refinement of MILP relaxations is inspired by similar algorithms from \cite{Gerdts2022} and \cite{Burlacu2019}. The idea is outlined in Algorithm \ref{alg:refine} and can be briefly described in the following way. First, some coarse discretizations $\mathcal{G}_i$ of the corresponding time intervals $I_{t_{i,min}}^{t_{i,max}},\ i \in \mathcal{N},$ are chosen. The algorithm then detects the subintervals where a predefined error bound is violated with the aim of refining them. In the next iteration step, the new refined time discretizations are examined. The whole procedure is repeated until the error bound is satisfied for all $i \in \mathcal{N}$ \cite{Gerdts2022}. As to the termination condition, it was proven in \cite{Burlacu2019} that the number of iteration steps is guaranteed to be finite provided that problem \eqref{prob:SP} is feasible, and all nonlinear functions are continuous.

\begin{algorithm}
\caption{Global optimization through adaptive refinement of MILP relaxations}
\label{alg:refine}
\begin{algorithmic}[1]
\Require An MINLP \eqref{prob:SP}
, initial tolerance $\epsilon_0$, the maximal linearization error $\epsilon$ and the number of new discretization points $L$.
\Ensure If \eqref{prob:SP}
is feasible, the algorithm returns an optimal solution $t^\ast_i,\ i \in \mathcal{N},$ of \eqref{prob:MILP}, which is an MILP relaxation of \eqref{prob:SP}, such that $\lvert D_i(t^\ast_i) - D_{i,pw} \rvert \leq \varepsilon$, $i \in \mathcal{N},$ and the corresponding objective value is minimal for any admissible point of \eqref{prob:SP}.

\State Choose initial discretizations $\mathcal{G}^0_i,\ i \in \mathcal{N},$ such that the corresponding nonlinearities do not exceed the initial tolerance value.
\State Set $r \gets 0$
\Repeat
    \State 
        Construct an MILP relaxation $\Pi^r$ of \eqref{prob:SP} from $\mathcal{G}^r_i,\ i \in \mathcal{N}$.
\State Solve $\Pi^r$.
\If{$\Pi^r$ is \textit{feasible}}
    \State Set $t_i^r,\ i \in \mathcal{N} \gets$ optimal solution of $\Pi^r$ 
    \State \multiline{
        Set $D^r_{i,pw},\ i \in \mathcal{N} \gets$ linear approximation of nonlinear function $D_i(t_i^r)$}
\Else
    \State \Return \textit{infeasible}
\EndIf    
    \State Set stop $\gets$ true
\ForAll{$i \in \mathcal{N}$} 
    \State Set $err^r \gets \lvert D_i(t^r_i) - D^r_{i,pw} \rvert$
\If{$err^r > \epsilon$}
    \State \multiline{%
        Set $n_0 \gets$ the index of the corresponding nonzero pointer $z_{i,n_0}^r$ from the solution of $\Pi^r$}        
    \State \multiline{%
        Set $h^r \gets \frac{t^r_{i,n_0+1} - t^r_{i,n_0}}{L+1}$} 
    \State \multiline{%
        Set $\mathcal{G}^r_{new} \gets \{t^r_{i,n_0} + h^r, t^r_{i,n_0} + 2h^r, \dots,  t^r_{i,n_0} + Lh^r\}$}         
    \State \multiline{%
        Set $\mathcal{G}_i^{r+1} \gets \mathcal{G}_i^{r} \cup \mathcal{G}^r_{new}$}    
    \State Set stop $\gets$ false 
\Else
    \State Set $\mathcal{G}_i^{r+1} \gets \mathcal{G}_i^{r} $ 
\EndIf
\EndFor
\State Set $r \gets r+1$
\Until stop 
\State \Return  $t^{r-1}_i,\ i \in \mathcal{N}$
\end{algorithmic}
\end{algorithm}

\subsection{Third Step: Computation of Optimal Trajectories}
\label{subsec:opt_traj}
In the last step, the optimal trajectories are computed, see  \cite{Bellman1957}. 
First, the optimal starting times $t^{\ast}_i, \ i \in \mathcal{N},$ are rounded up to the next node in the time grid. This step is not necessary, but is performed here due to safety reasons in order to avoid motions' overlaps. Without performing this step, interpolation of the value function in time is required in addition. We denote those new starting times by $t_{k_i}, \ i \in \mathcal{N}$.
For each $\ i \in \mathcal{N}$, optimal trajectories can be computed as follows:
\begin{align*}
	u^{\ast}_i(t_k) &= \underset{u \in U}{\arg \min} \{h + V_h(t_{k+1}, f(x^{\ast}_i(t_k),u))\},\\
	x^{\ast}_i(t_{k+1}) &= f(x^{\ast}_i(t_k),u^{\ast}_i(t_k)),
	\quad k = k_i,\ k_i + 1, \ k_i + 2,... \ .
\end{align*}
If $V_h(t_{k+1}, f(x^{\ast}_i(t_k),u)) = 0$ for some $k$, $	x^{\ast}_i(t_{k+1})$ is the last point in the trajectory meaning that the algorithm terminates.

Note that error estimates are possible, and that the error in $V_h$ can be controlled.

\section{Numerical Results}
\label{sec:numerics}

All numerical experiments were conducted on an Intel(R) Core i7 processor with a 2.80GHz-CPU and an 8GB-RAM using the robust and efficient GUROBI MILP-solver \cite{Gurobi} based on the Branch-and-Bound algorithm. For reasons of space, the idea behind this algorithm is not explained here. The interested reader is instead referred to \cite{Winston2004}.

Our approach was extensively tested for up to $N = 8$ VTOLs to show its feasibility. A periodic problem was considered. The time-space domain $\Omega$ with $n=2$ and the control set $U$ were replaced by their finite discretizations. The value function $V_h$ was then obtained by its pointwise evaluation at the corresponding discretization points. All results presented below relate to the problem setting described as follows.

In the state space $\Omega = \interval{-200}{200} \times \interval{-200}{200} \subset \mathbb{R}^2$ with $\Omega_T$ lying in its center, we consider $N=8$ VTOLs and $P=4$ dynamic obstacles.
VTOLs in our consideration have the following initial positions:
\begin{align*}
\overline{x}_{1} &= \overline{x}_{5} = \{-192, - 192\}\si{m}, \\
\overline{x}_{2} &= \overline{x}_{6} = \{192, - 192\}\si{m}, \\
\overline{x}_{3} &= \overline{x}_{7} = \{192,  192\}\si{m}, \\
\overline{x}_{4} &= \overline{x}_{8} = \{-192, 192\}\si{m}.
\end{align*}
They move with a constant velocity $v = 10 \si{m}/\si{s}$. The circular obstacles, with $R_O = 64\si{m}$ radius each, move periodically around the target region $\Omega_T$ along a circular trajectory with radius $R_C = 100\si{m}$. The initial positions of their centers are given by $\overline{x}_1 = [R_C, 0]^T, \overline{x}_2 = [0, R_C]^T, \overline{x}_3 = [-R_C, 0]^T$ and $\overline{x}_4 = [0, -R_C]^T$, respectively. Having a period of $80\si{s}$ (full rotation around the target), their velocity, which is assumed to be constant for convenience, is approximately equal to $7.85\si{m}/\si{s}$.

\begin{figure}[htbp!]
\begin{center}
\includegraphics[scale=0.8]{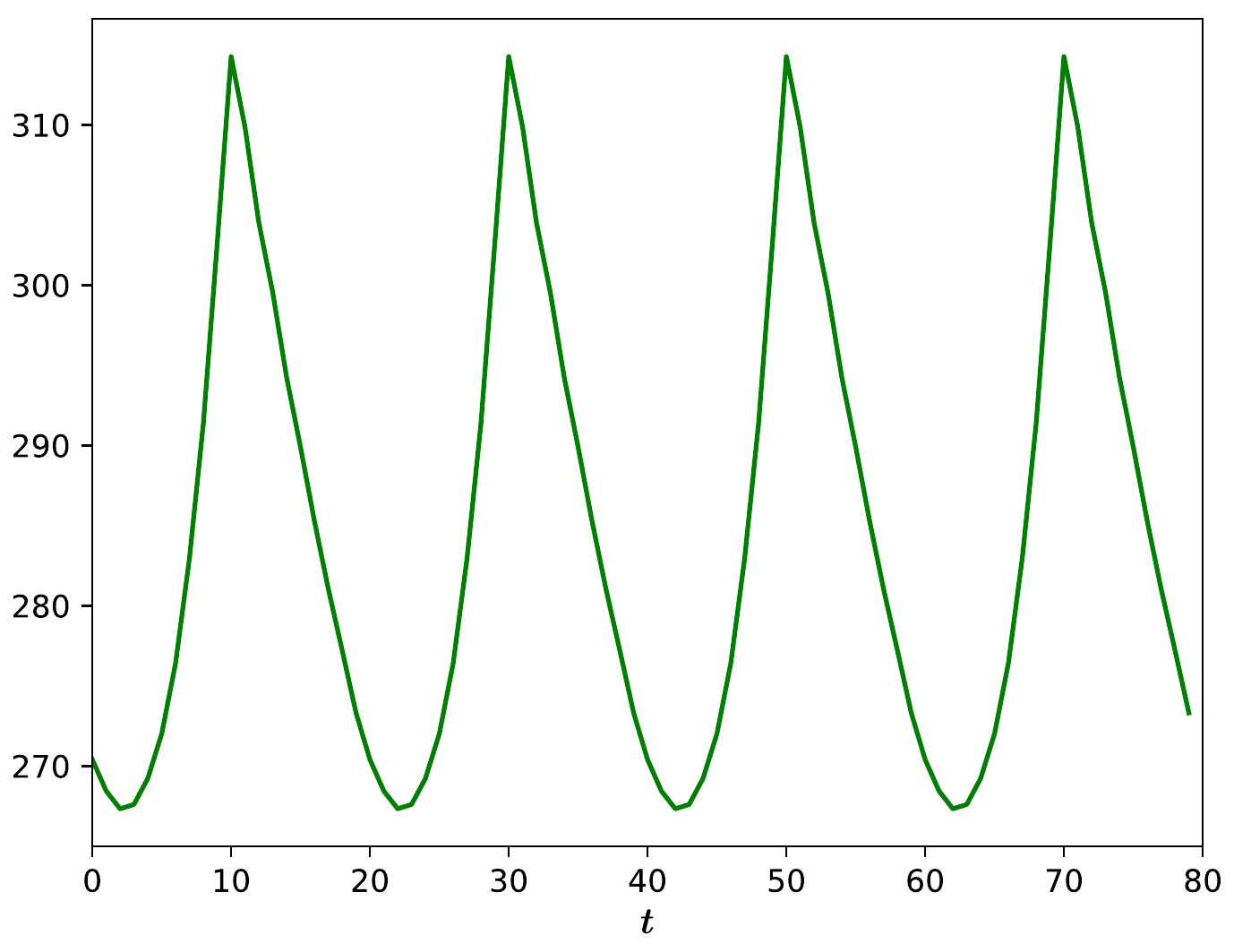}
\end{center}
\caption{Shortest-path function of a VTOL starting at point $\overline{x} = [192,192]^T \si{m}$}
\label{fig:val_fun_per}
\end{figure}

Note that the value function $V_h(t,x)$ evaluated at some point $\overline{x}$ corresponds to a time-dependent function, which indicates the length of the shortest path from that point to the target region $\Omega_T$. Figure \ref{fig:val_fun_per} illustrates this function for a VTOL starting at $\overline{x} = [ 192, 192 ]^T \si{m}$. Since all four obstacles move periodically with a constant and identical velocity, it suffices to analyze the plot during the first 20\si{s}. The shortest path is observed when the VTOL starts moving after  $2\si{s}$, while the longest path has a length of nearly 315\si{m} and refers to the case where VTOL commences its motion at time $t = 11\si{s}$. This observation can be explained by the fact that the VTOL does not collide with any of the obstacles if it starts at $t= 2\si{s}$ and moves straight towards the target. On the contrary, if the VTOL starts at $t= 11 \si{s}$, it has to avoid a collision with an obstacle and therefore flies around it. For that reason, the path is in the second case much longer than in the first one. After 20\si{s}, the next obstacle takes place of the previous one, and the same arguments apply. 

\begin{figure}[htbp!]
\begin{center}
\includegraphics[scale=0.8]{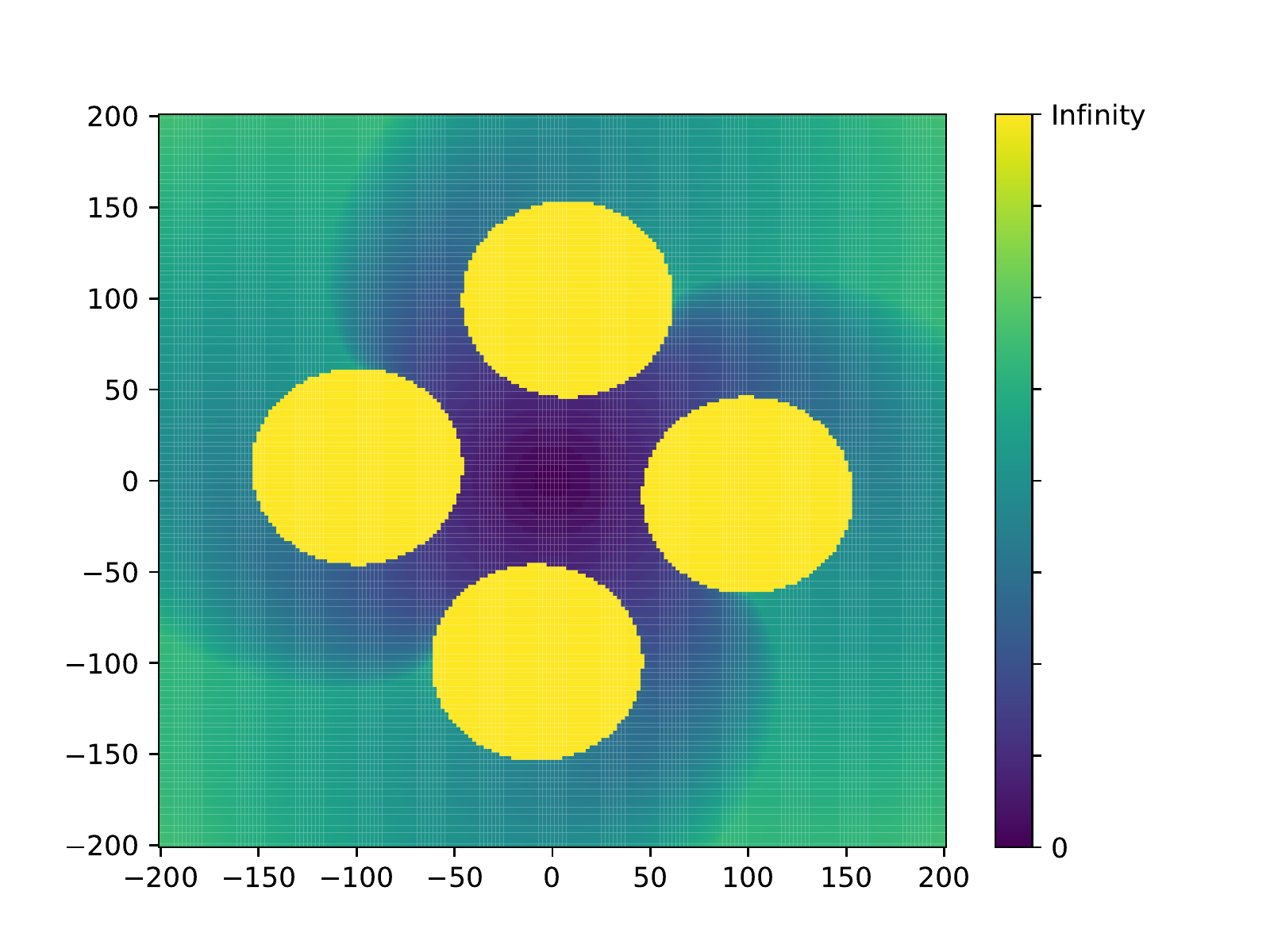}
\end{center}
\caption{Value function $V_h(t,x)$ evaluated at time $t_0=0$}
\label{fig:val_fun_x}
\end{figure}

Next, we present a color map shown in Figure \ref{fig:val_fun_x}, which visualizes the value function evaluated at the initial time point $t_0=0$. To be more specific, for any point $x \in \Omega(t_0) = \interval{-200}{200} \times \interval{-200}{200}$ this plot indicates the length of the shortest path (taking into account obstacle avoidance) from a VTOL located at this point to the target region $\Omega_T$ with radius $R_T = 10 \si{m}$. Obviously, if the initial position of a VTOL lies somewhere within any of the obstacles, the value function $V_h$ is equal to an infinite value according to our convention. Therefore, the four circles in Figure \ref{fig:val_fun_x} are colored yellow. Similarly, if the VTOL starts at a point lying within the target set, its distance to this set is zero. This explains why $\Omega_T$ illustrated as a circle in the middle of the color map has the darkest color. The last interesting observation arises from the semicircular regions depicted in Figure \ref{fig:val_fun_x}, which take such a shape because of the obstacles' round shape. If the initial position of a VTOL lies somewhere within such blue region, the path from the VTOL to $\Omega_T$ would be rather short due to the possibility of a collision-free flight along a nearly straight trajectory.

In order to determine the average computational time required by the GUROBI MILP solver to solve problem \eqref{prob:MILP}, 10 experiments were run with the same parameter setting. For $\alpha = 1$, the solver needed on average 19.34\si{s} to find a numerical solution of the final MILP \ref{prob:MILP}. In case that $\alpha = 10$, the computational time reduced drastically to 3.65\si{s}. The achieved results yield that the solver runtime decreases with a growing $\alpha$. This relation is caused by the fact that larger values of the penalization term yield smaller search trees in the Branch-and-Bound algorithm and hence shorter computational time.

Optimal schedules obtained by the solver for $\alpha = 1$ and $\alpha = 10$ are shown in Tables \ref{table:alpha_1} and \ref{table:alpha_10}, respectively. A closer look at the last column of Tables \ref{table:alpha_1} and \ref{table:alpha_10} reveals that the larger the $\alpha$ values are, the smaller the flight duration of a single VTOL becomes. This observation is totally plausible, since the VTOLs are penalized more in the second case and therefore try to complete their motions as quickly as possible. While the duration of a \textit{single} mission decreases with a growing $\alpha$, the \textit{total} task duration shows an opposite tendency. In our problem setting, the last VTOL reached the target set after $225.57\si{s}$ for $\alpha=1$, whereas the overall mission was completed only after $247.885\si{s}$ for $\alpha=10$. This phenomenon can be explained by the fact that each VTOL generally tends to wait longer when $\alpha$ increases in order to reduce its own flight duration. As a result, the overall mission time increases.

\renewcommand{\arraystretch}{1.2}
\begin{table}[ht!]
	\centering
\begin{tabular}{|c|c|c|c|c|}
 \hline
  Starting Number & Original Number & Starting Time [\si{s}] & End Time [\si{s}] & Duration [\si{s}] \\
 	\hline
 	\hline
 	1 & 4 & 0 & 27.0094 & 27.0094 \\
	\hline
	2 & 3 & 27.0094 & 55.1403 & 28.131 \\
	\hline
	3 & 1 & 57.1793 & 84.9432 & 27.7639 \\
	\hline
	4 & 2 & 84.9432 & 111.925 & 26.9819 \\
	\hline
	5 & 5 & 111.925 & 142.105 & 30.1801 \\
	\hline
	6 & 8 & 142.105 & 168.331 & 26.2255 \\
	\hline
	7 & 6 & 168.331 & 197.768 & 29.4372 \\
	\hline
	8 & 7 & 197.768 & 225.57 & 27.8021 \\
	\hline
\end{tabular}
\caption{Optimal Schedule for $N=8$ VTOLs with $\alpha=1$}
\label{table:alpha_1}
\end{table}

\begin{table}[ht!]
	\centering
\begin{tabular}{|c|c|c|c|c|}
 \hline
  Starting Number & Original Number & Starting Time [\si{s}] & End Time [\si{s}] & Duration [\si{s}] \\
 	\hline
 	 \hline
 	1 & 1 & 2.63158 & 29.2895 & 26.658 \\
	\hline
	2 & 3 & 39.4737 & 66.4896 & 27.0159 \\
	\hline
	3 & 4 & 66.4896 & 94.339 & 27.8494 \\
	\hline
	4 & 2 & 98.9736 & 126.316 & 27.3422 \\
	\hline
	5 & 8 & 126.316 & 153.973 & 27.6573 \\
	\hline
	6 & 6 & 158.427 & 185.965 & 27.5377 \\
	\hline
	7 & 5 & 185.965 & 213.464 & 27.4995 \\
	\hline
	8 & 7 & 221.053 & 247.855 & 26.8026 \\
	\hline
\end{tabular}
\caption{Optimal Schedule for $N=8$ VTOLs with $\alpha=10$}
\label{table:alpha_10}
\end{table}

Finally, a 2D visualization was created in order to illustrate the results from above for two different values of the penalization term, namely for $\alpha=1$ and $\alpha=10$. The respective videos can be found under:

\url{https://youtu.be/QgMEM1QQP9M} \\
\indent
\url{https://youtu.be/BQ39fRMV-ew}

\section{Conclusion}
\label{sec:conclusion}
The underlying work presented an efficient numerical framework for scheduling of multiple autonomous VTOLs and their simultaneous optimal trajectory planning. Formulated as a bilevel problem, the model is recast into a single-level one by means of a value function approach. This technique involves such advanced concepts as the Hamilton-Jacobi-Bellman equation originating from dynamic programming theory and the Kru\u{z}kov transformation. The resulting MINLP was piecewise linearized with the help of an elegant technique allowing to track the linearization error. The solution of the final MILP was then determined numerically by the GUROBI solver based on the Branch-and-Bound algorithm. The numerical experiments yield the feasibility of the framework presented in this article.

Possible model extensions encompass:
\begin{itemize}
 \item Extension of the optimal control problem at the lower level to cope with multi-phase systems.
 \item Consideration of alternative scheduling mixed-integer problems at the upper level as, for example, vehicle routing.
 \item Integration of other objectives at the upper level with the aim of maximizing the mission success, minimizing the fuel consumption etc.
 \item Integration of coupling effects by introducing common constraints and/or common objectives.
 \item Consideration of multi-stage problems with early starts.
\end{itemize}

\section*{Funding}
This material is based upon work supported, inter alia, by dtec.bw–Digitalization and Technology Research Center of the Bundeswehr [MissionLab] and by the Air Force Office of Scientific Research, Air Force Materiel Command, USAF, under Award No. FA8655-20-1-7026. Any opinions, findings, and conclusions or recommendations expressed in this publication are those of the author(s) and do not necessarily reflect the views of the funders.

\bibliographystyle{tfnlm}


\end{document}